\documentclass[final,1p,times]{elsarticle}

\usepackage{amsmath, amsthm}
\usepackage{amssymb}
\usepackage{hyperref}
\usepackage{esint}
\usepackage{color}
\usepackage{mathrsfs}
\usepackage{bbold}
\usepackage{enumerate}

\makeatletter
\providecommand*{\dif}%
   {\@ifnextchar^{\DIfF}{\DIfF^{}}}
\def\DIfF^#1{%
   \mathop{\mathrm{\mathstrut d}}%
      \nolimits^{#1}\gobblespace
}
\def\gobblespace{%
   \futurelet\diffarg\opspace}
\def\opspace{%
   \let\DiffSpace\!%
   \ifx\diffarg(%
      \let\DiffSpace\relax
   \else
      \ifx\diffarg\[%
         \let\DiffSpace\relax
      \else
         \ifx\diffarg\{%
            \let\DiffSpace\relax
         \fi\fi\fi\DiffSpace}

\theoremstyle{plain}
\newtheorem{theorem}{Theorem}[section]

\newtheorem{lemma}[theorem]{Lemma}

\newtheorem{corollary}[theorem]{Corollary}

\theoremstyle{definition}
\newtheorem{definition}[theorem]{Definition}

\theoremstyle{remark}
\newtheorem{remark}[theorem]{Remark}

 \usepackage{amsthm}

\begin{document}
\begin{frontmatter}



\title{Two-weight Norm Estimates for Multilinear Fractional Integrals in Classical Lebesgue Spaces$^*$}

\author[add1,add2]{Vakhtang Kokilashvili}
\ead{kokil@rmi.ge}

\address[add1]{Department of Mathematical Analysis, A. Razmadze Mathematical Institute, I. Javakhishvili Tbilisi State University,
2. University Str., 0186 Tbilisi, Georgia}
\address[add2]{International Black Sea University, 3 Agmashenebeli Ave., Tbilisi 0131, Georgia }

\author[add5]{Mieczyslaw Mastylo}
\address[add5]{Faculty of Mathematics and Computer Science, Adam Mickiewicz University Pozna\'n,
and Institute of Mathematics,  Polish Academy of Sciences (Pozna\'n branch),  Umultowska 87, 61-614 Pozna\'n, Poland}
\ead{mastylo@amu.edu.pl}

\author[add1,add3,add4]{Alexander Meskhi\corref{cor}}
\address[add3]{Department of Mathematics,  Faculty of Informatics and Control
Systems,  Georgian Technical University, 77, Kostava St., Tbilisi, Georgia}
\ead{alex72meskhi@yahoo.com; meskhi@rmi.ge}
\cortext[cor]{This paper was publish in  {\em Preprints of Abdus Salam School of Mathematical Sciences}, No. 481, at October 29, 2012 (see http://sms.edu.pk/journals/preprint/pre\_481.pdf)}

\begin{abstract}
We derive criteria governing two-weight estimates  for multilinear fractional integrals  and appropriate maximal functions. The two and one weight problems for multi(sub)linear strong fractional maximal operators are  also studied; in particular,  we derive necessary and sufficient conditions guaranteeing the trace type inequality for this operator. We also establish the Fefferman-Stein type inequality, and obtain
one-weight criteria when a weight function is of product type. As a consequence, appropriate results for multilinear Riesz potential operator with product kernels follow.

\end{abstract}

\begin{keyword}
Multilinear fractional integrals; multi(sub)linear strong fractional maximal operators; boundedness; two-weight inequality; trace inequality; Fefferman--Stein inequality.

\MSC[2010] Primary  42B25 \sep Secondary \sep 42B35

\end{keyword}

\end{frontmatter}



\section{ Introduction}

Necessary and sufficient conditions governing the two-weight inequality for multilinear fractional integral operators are established. Our results involve, for example, criteria ensuring the trace type inequality and one-weight inequality for multi(sub)linear fractional maximal operators. Fefferman--Stein type inequality is also established for these operators. It should be also emphasized that appropriate results for multilinear Riesz potential operators follow.

Historically, multilinear fractional integrals were introduced in the papers by L. Grafakos \cite{Graf1},  C. Kenig and E. Stein \cite{KeSt},  L. Grafakos and N. Kalton \cite{GrKa}. In particular, these works deal with the operator

$$ B_{\alpha}(f,g)(x) =\int\limits_{{\Bbb{R}}^n} \frac{f(x+t)g(x-t)}{|t|^{n-\alpha}} dt, \;\;\; 0< \alpha<n. $$

In the mentioned  papers it was proved that if $\frac{1}{q} = \frac{1}{p}-\frac{\alpha}{n}$, where $\frac{1}{p}=\frac{1}{p_1}+\frac{1}{p_2}$, then $B_{\alpha}$ is bounded from $L^{p_1}\times L^{p_2}$  to  $L^q$.

As a tool to understand $B_{\alpha}$, the operators
$$ {\mathcal{I}}_{\alpha}(\overrightarrow{f}) (x) = \int\limits_{({\Bbb{R}}^n)^m}\frac{f_1(y_1)\cdots  f_m(y_m)}{(|x-y_1|+ \cdots + |x-y_m|)^{mn-\alpha}} d\overrightarrow{y}, $$
where $x\in {\Bbb{R}}^n$, $0<\alpha<nm$, $\overrightarrow{f}:= (f_1, \cdots, f_m)$, $\overrightarrow{y}:=(y_1, \cdots, y_m)$,  were studied as well. Corresponding maximal operator is given by (see \cite{Mo})
$${\mathcal{M}}_{\alpha} (\overrightarrow{f})(x) =\sup_{Q\ni x} \prod_{i=1}^m \frac{1}{ |Q|^{1-\alpha/(nm)}} \int\limits_{Q} |f_i(y_i)| d y_i, \;\;\; 0\leq \alpha< mn,  $$
where $|Q|$ denotes the volume of the cube $Q$ with sides parallel
to  the coordinate axes.

This operator for $\alpha=0$ was introduced and studied in \cite{LOPTT-G}.

It  can be immediately checked that
$$ {\mathcal{I}}_{\alpha} (\overrightarrow{f})(x) \geq c_n {\mathcal{M}}_{\alpha}(\overrightarrow{f}) (x) , \;\;\; f_i\geq 0, \;\;\; i=1, \cdots , m. $$

In the sequel the following notation will be used:

$$   \overrightarrow{p}:=(p_1, \cdots, p_m), \;\;\; \overrightarrow{w}= (w_1, \cdots, w_m),$$
where $p_i$ are constants ($0< p_i<\infty$) and $w_i$ are a.e. positive functions defined on Euclidean space. It will be also assumed that
$$ \frac{1}{p}= \sum_{i=1}^m \frac{1}{p_{i}}. $$

\vskip+0.1cm

Throughout the paper we use the notation ${\mathcal{Q}}$ to denote the family of all cubes in $\mathbb{R}^n$ with sides parallel to the coordinate axes.

An almost everywhere positive and locally integrable function defined on ${\Bbb{R}}^n$ is called a weight.

\vskip+0.2cm

Let $0<r<\infty$ and let $\rho$ be a weight on ${\Bbb{R}}^n$. We denote by $L^{r}_{\rho}({\Bbb{R}}^n)$ a class of all measurable functions $f$ on ${\Bbb{R}}^n$ such that
$$ \|f\|_{L^r_{\rho}({\Bbb{R}}^n)}:= \bigg( \int\limits_{{\Bbb{R}}^n} |f(x)|^r \rho(x) dx \bigg)^{1/r} < \infty. $$

\vskip+0.2cm
\begin{definition} [vector Muckenhoupt condition]  Let $1\leq p_i <\infty$ for $i= 1, \cdots, m$. Let $w_i$ be weights on ${\Bbb{R}}^n$, $i=1, \cdots, m$. We say that  $\overrightarrow{w}\in A_{\overrightarrow{p}}({\Bbb{R}}^n)$ (or simply $\overrightarrow{w}\in A_{\overrightarrow{p}}$) if
$$ \sup_{Q\in {\mathcal{Q}}} \bigg( \frac{1}{|Q|} \int\limits_{Q} \prod_{i=1}^m w_i^{p/p_i}\bigg)^{1/p} \prod_{i=1}^m \bigg( \frac{1}{|Q|}\int\limits_{Q} w_i^{1-p_i'}\bigg)^{1/p_i'}< \infty. $$
When $p_i=1$, $\Big( \frac{1}{|Q}\int\limits_{Q}w_i^{1-p_i'}\Big)^{1/p_i'}$ is understood as $\big( \inf_Q w_i\big)^{-1}$.
\end{definition}

\begin{remark} In the linear case ($m=1$) the class  $A_{\overrightarrow{p}}$ coincides with the well- known Muckenhoupt class $A_p$.
\end{remark}

"Let us recall that,  by definition, $A_{\infty}= \cup_{i=1}^{\infty}A_p$.


Further, let us recall that, by definition, $A_{\infty}= \cup_{i=1}^{\infty}A_p$.

\vskip+0.2cm

\begin{definition} [vector Muckenhoupt-Wheeden condition]  Let $1\leq p_i <\infty$ for $i= 1, \cdots, m$. Suppose that $p<q<\infty$.  We say that $\overrightarrow{w}= (w_1, \cdots, w_m)$ satisfies $A_{\overrightarrow{p},q}({\Bbb{R}}^n)$ condition $(\overrightarrow{w}\in A_{\overrightarrow{p},q})$ if
$$ \sup_{Q} \bigg( \frac{1}{|Q|} \int\limits_{Q} \Big( \prod_{i=1}^m w_i\Big)^q \bigg)^{1/q} \prod_{i=1}^m \bigg( \frac{1}{|Q|}\int\limits_{Q} w_i^{-p_i'}\bigg)^{1/p_i'}< \infty. $$
\end{definition}

{\bf Theorem A (\cite{LOPTT-G}). } {\em Let $1<p_i<\infty$, $i=1, \cdots, m$. Suppose that  $w_i$ are weights on ${\Bbb{R}}^n$. Then the operator ${\mathcal{M}}$ is bounded from $L^{p_1}_{w_1}({\Bbb{R}}^n) \times \cdots \times L^{p_m}_{w_m}({\Bbb{R}}^n)$ to
$L^{p}_{\prod_{i=1}^m w^{p/p_i}_i}({\Bbb{R}}^n)$ if and only if $\overrightarrow{w}\in A_{\overrightarrow{p}}({\Bbb{R}}^n)$.}

\vskip+0.1cm

{\bf Theorem B.} (\cite{Mo}) {\em Let $1<p_1, \cdots, p_m<\infty$, $0<\alpha <mn$, $\frac{1}{m} < p<\frac{n}{\alpha}$. Suppose that $q$ is an exponent satisfying the condition $\frac{1}{q}= \frac{1}{p}-\frac{\alpha}{n}$. Suppose that  $w_i$ are a.e. positive functions  on ${\Bbb{R}}^n$ such that $w_i^{p_i}$ are weights. Then the inequality
$$\bigg( \int\limits_{{\Bbb{R}}^n} \Big( \big| {\mathcal{N}}_{\alpha} (\overrightarrow{f})  (x) \big| \prod_{i=1}^m w_i(x) \Big)^q    dx \bigg)^{1/q} \leq C \ \prod_{i=1}^m \bigg( \int\limits_{{\Bbb{R}}^n} \big( |f_i(y_i)|w_i\big)^{p_i} dy_i \bigg)^{1/p_i}, $$
holds, where ${\mathcal{N}}_{\alpha}$ is ${\mathcal{I}}_{\alpha}$ or ${\mathcal{M}}_{\alpha}$, if and if $\overrightarrow{w} \in A_{\overrightarrow{p}, q}({\Bbb{R}}^n)$. }
\vskip+0.2cm

Regarding the two-weight result it is known the following result:
\vskip+0.2cm
{\bf Theorem C.} {\em Let $1<p_1, \cdots, p_m<\infty$ and let  $0<\alpha <mn$. Let $w_i$ and $u$ be a.e. positive functions on ${\Bbb{R}}^n$.  Further, assumed that $q$ is an exponent satisfying the condition $\frac{1}{m}< p\leq q<\infty$. Suppose also that

\rm{(i)} $q>1$ and a pair $(u, \overrightarrow{w})$ satisfies the condition

\begin{equation}\label{eq:powerbump1}
\sup_{Q}|Q|^{\alpha/n+ 1/q-1/p} \bigg( \frac{1}{|Q|}\int\limits_{Q} u^{qr}(x) dx \bigg)^{1/qr} \prod_{i=1}^m \bigg( \frac{1}{|Q|}\int\limits_{Q} w_i^{-p_i' r}(y_i) dy_i\bigg)^{1/p_i' r} <\infty
\end{equation}
for some $r>1$;

\rm{(ii)} $q\leq 1$ and $(u, \overrightarrow{w})$ satisfies the condition

\begin{equation}\label{eq:powerbump2}
\sup_{Q}|Q|^{\alpha/n+ 1/q-1/p} \bigg( \frac{1}{|Q|}\int\limits_{Q} u^{q}(x) dx \bigg)^{1/qr} \prod_{i=1}^m \bigg( \frac{1}{|Q|}\int\limits_{Q} w_i^{-p_i' r}(y_i) dy_i\bigg)^{1/p_i' r} <\infty
\end{equation}
for some $r>1$. Then the inequality
$$ \bigg( \int\limits_{ {\Bbb{R}}^n } \big( |{\mathcal{I}}_{\alpha} (\overrightarrow{f})(x)| u(x) \big)^q dx\bigg)^{1/q}
\leq c \prod_{i=1}^m \bigg( \int\limits_{{\Bbb{R}}^n} \big( |f_i(x)|w_i(x)\big)^{p_i} dx \bigg)^{1/p_i} $$
holds for all $\overrightarrow{f} \in L^{p_1}_{w_1^{p_1}}({\Bbb{R}}^n) \times \cdots \times L^{p_m}_{w_m^{p_m}}({\Bbb{R}}^n)$.}

\vskip+0.2cm

\begin{remark}
The two-weight problem for sublinear fractional maximal operator
$$ \big( M_{\alpha}f\big) (x) =\sup_{Q\ni x} \frac{1}{|Q|^{1-\alpha/n}} \int\limits_{Q} f(y) dy, $$
has been already solved. We mention the papers by E. Sawyer \cite{Saw1} for the conditions  involving the operator itself and  by R. L. Wheeden  \cite{Wh} for the Gabidzashvili-Kokilashvili type conditions (see also   \cite{KokMim}, \cite{KoKr} for the latter condition). For  the solution of the two-weight problem for the Riesz fractional integral operator
$$(I_{\alpha} f)(x) = \int\limits_{{\Bbb{R}}^n}\frac{f(y)}{|x-y|^{n-\alpha}} dy  $$
under different conditions we refer to the paper \cite{SawPot} and the  monograph  \cite{KoKr}.
\end{remark}
\begin{remark}
Conditions \eqref{eq:powerbump1} and \eqref{eq:powerbump2} are known as the ''power bump'' conditions (see \cite{Pe}) for the linear case).
\end{remark}

\vskip+0.2cm

\begin{remark}
Some ''power bump'' sufficient conditions governing the two-weight inequality for ${\mathcal{I}}_{\alpha}$ were derived in \cite{ShTa}.
\end{remark}

\vskip+0.1cm
{\bf Corollary A.} {\em Let $0<\alpha <mn$ and let  $1<p_i<\infty$, $i=1, \cdots, m$. Assume that  $\frac{1}{m}<p\leq q <\infty$. Suppose that  $u, w_1, \cdots w_m$ are a.e. positive functions on ${\Bbb{R}}^n$ such that
$u^q, w_1^{-p_1'}, \cdots, w_m^{-p'_m} \in A_{\infty} ({\Bbb{R}}^n)$. Then the inequality

\begin{equation}\label{eq:two-weightInequality}
\bigg( \int\limits_{{\Bbb{R}}^n} \Big( |{\mathcal{N}}_{\alpha}(\overrightarrow{f})(x)| u(x)\Big)^q dx  \bigg)^{1/q}
\leq C \prod_{i=1}^m \bigg( \int\limits_{{\Bbb{R}}^n} \big( |f_i(y_i)|w_i(y_i)\big)^{p_i} dy_i \bigg)^{1/p_i},
\end{equation}
where ${\mathcal{N}}_{\alpha}$ is ${\mathcal{I}}_{\alpha}$ or ${\mathcal{M}}_{\alpha}$ holds if and only if

\begin{equation} \label{eq:two-weightCondition}
\sup_{Q\in {\mathcal{Q}}}  |Q|^{\alpha/n+1/q-1/p} \bigg( \frac{1}{|Q|} \int\limits_{Q} u^{q}(x) dx \bigg)^{1/q} \prod_{i=1}^m \bigg( \frac{1}{|Q|} \int\limits_{Q} w_i^{-p'_i}(y_i) dy_i \bigg)^{1/p'_i}< \infty.
\end{equation}}

\vskip+0.2cm

{\bf Theorem D.} (\cite{Mo}) {\em Let $0\leq \alpha < mn$ and let $1\leq p_i < \infty$, $1= 1, \cdots, m$. Suppose that $\frac{1}{m} < p \leq q <\infty$. Then the weak type inequality
$$ \|u  {\mathcal{M}}_{\alpha}(\overrightarrow{f}) \|_{L^{q, \infty}({\Bbb{R}}^n)} \leq c \prod_{i=1}^{m} \| w_i f_i\|_{L^{p_i}({\Bbb{R}}^n)}$$
holds if and only if \eqref{eq:two-weightCondition} satisfied.}

\vskip+0.2cm

The next statement shows that for ${\mathcal{N}}_{\alpha}= {\mathcal{M}}_{\alpha}$, only the ''power bump'' condition with respect to the right-hand side weight is sufficient for the strong--type inequality for ${\mathcal{M}}_{\alpha}$.
\vskip+0.2cm

{\bf Theorem E.} {\em Let $0\leq \alpha< mn$, $1<p_1, \cdots, p_n$, $\frac{1}{m}< p \leq q < \infty$.  If $(u, \overrightarrow{w})$ satisfies \eqref{eq:powerbump2}, then inequality
\eqref{eq:two-weightInequality} holds for ${\mathcal{N}}_{\alpha}= {\mathcal{M}}_{\alpha}$.}
\vskip+0.2cm

{\bf Corollary B.} {\em Let $\alpha$, $p_i$, $q$ satisfy the conditions of   Theorem E. Suppose that $w_i^{-p_i'} \in A_{\infty}$. Then the inequality \eqref{eq:two-weightInequality} holds for ${\mathcal{N}}_{\alpha}= {\mathcal{M}}_{\alpha}$ if and only if \eqref{eq:two-weightCondition} holds.}
\vskip+0.2cm

\vskip+0.2cm

In  \cite{KoMaMe}, based on Corollary B,  the authors of this paper characterized the inequality \eqref{eq:two-weightInequality} for ${\mathcal{N}}_{\alpha}={\mathcal{I}}_{\alpha}$ and $w_i \equiv \; const$ without any additional restriction on $u$ under the D. Adams \cite{Ad} type condition. In particular, they proved the  next statement.
\vskip+0.2cm

{\bf Theorem F.}  {\em Let $1<p_i<\infty$, $i=1, \cdots, m$. Assume that  $\alpha<n/p$ and $p< q <\infty$. Then the following conditions are equivalent:

\rm{(i)} inequality \eqref{eq:two-weightInequality}  holds for ${\mathcal{N}}_{\alpha}={\mathcal{I}}_{\alpha}$ and $w_i \equiv \; const$;

\rm{(ii)} the  inequality

$$ \|u  {\mathcal{M}}_{\alpha}(\overrightarrow{f}) \|_{L^{q, \infty}({\Bbb{R}}^n)} \leq c \prod_{i=1}^{m} \|  f_i\|_{L^{p_i}({\Bbb{R}}^n)}$$
holds;

\rm{(iii)} the condition

\begin{equation*}
\sup_{Q\in {\mathcal{Q}}}  \bigg(\int\limits_{Q} u^q(x) dx \bigg)^{1/q} |Q|^{\alpha-n/p}<   \infty
\end{equation*}
 is satisfied.}
\vskip+0.2cm

\begin{remark}
Finally we mention that Fefferman-Stein type inequality for the multi(sub)linear operator ${\mathcal{M}}_{\alpha}$ was studied in \cite{Pr}.
\end{remark}

\section{Preliminaries}

Let ${\mathcal{D}}$ be the set of all dyadic cubes in ${\Bbb{R}}^n$. It is known that each $Q\in {\mathcal{D}}$ is the
union of $2^n$ non-overlapping dyadic cubes (for the definition and some properties of the
dyadic intervals see, for instance, \cite{GCRF}, P. 136).

\begin{definition} We say that a weight function $\rho$ satisfies the dyadic reverse doubling condition $(\rho\in
RD^{(d)}({\Bbb{R}}))$ if there exists a constant $d>1$ such that

\begin{equation} \label{eq:reversedoubling}
d \rho(Q') \leq \rho(Q),
\end{equation}
for all $Q',Q\in {\mathcal{D}}$, where $Q'\subset Q$ and $|Q|=2^n |Q'|$.

Further, we say that a weight function $\rho$ satisfies the reverse doubling condition $(\rho\in RD({\Bbb{R}}^n))$
if  \eqref{eq:reversedoubling} holds for all cubes $Q',Q$, where $Q'$ is a subcube of $Q$ arising  dividing $Q$ by $2^n$ equal parts.
\end{definition}

We shall also need the next Carleson-H\"ormander  type embedding theorem regarding the dyadic intervals (see e.g., \cite{SaWh}, \cite{Ta}):
\vskip+0.2cm

\textsc{Theorem G.} {\em Let $1<r<q<\infty$ and let $\rho$ be a
weight function on ${\Bbb{R}}^n$ such that $\rho^{1-r'}$ satisfies the dyadic reverse doubling condition. Let $\{c_Q\}$ be non-negative numbers corresponding to dyadic cubes  $Q$ in ${\Bbb{R}}$. Then the following two statements are equivalent:

\rm{(i)} {\em There is a positive constant $C$ such that}
$$ \sum_{Q\in {\mathcal{D}}} c_Q\left({1\over |Q|}\int_Q g(x)dx\right)^q\le C \left(\int_{{\Bbb{R}}^n}g(x)^r \rho(x)dx\right)^{q/r} $$ {\em for all
non-negative} $g\in L^r_{\rho}({\Bbb{R}}^n)$;

\rm{(ii)} {\em There is a positive constant $C_1$ such that}
$$ c_Q \le C_1 |Q|^q\left(\int_Q \rho(x)^{1-r'}dx\right)^{-q/r'}$$
{\em for all  $Q\in {\mathcal{D}}$.}}

\vskip+0.1cm

This result yields the following corollary:

\vskip+0.1cm

{\bf Corollary C.} {\em Let $1<r<q<\infty$ and let $\rho$ be a weight function on ${\Bbb{R}}^n$ such that $\rho^{1-r'}$ satisfies the dyadic reverse doubling condition. Then the Carleson-H\"ormander type inequality
$$ \sum_{Q\in {\mathcal{D}}} \bigg( \int_{Q} \rho^{1-r'}(x)dx\bigg)^{-q/r'} \bigg(\int_{Q} f(x) dx
\bigg)^q\leq c \bigg(\int_{{\Bbb{R}}^n} f^r(x)\rho(x)dx\bigg)^{q/r}$$
holds for all non-negative $f\in L^r_{\rho}({\Bbb{R}}^n)$.}

\vskip+0.2cm

Let us recall some properties of vector Muckenhoupt condition:
\vskip+0.2cm

\begin{lemma}[\cite{LOPTT-G}, \cite{Mo}]\label{MuckenhouptProperties}
Let $1<p_i<\infty$, $i=1, \cdots, m$. Then

\rm{(i)} $$ \prod_{i=1}^m A_{p_i}({\Bbb{R}}^n) \subset A_{\overrightarrow{p}}({\Bbb{R}}^n);$$

\rm{(ii)}  if $\overrightarrow{w} \in A_{\overrightarrow{p},q}({\Bbb{R}}^n)$, where $1/q=1/p-\alpha/n$, $0<\alpha<mn$, then

$$ \bigg( \prod_{i=1}^m w_i \bigg)^q \in A_{mq}({\Bbb{R}}^n) \;\; \text{and} \;\; w_i^{-p'_i}\in A_{mp_i'}({\Bbb{R}}^n). $$
\end{lemma}

\section{Main Results}

Our aim in this paper is to improve the known results (see the statements above) regarding the two--weight strong-type inequality for ${\mathcal{N}}_{\alpha}$, where
${\mathcal{N}}_{\alpha}$ is  ${\mathcal{M}}_{\alpha}$ or ${\mathcal{I}}_{\alpha}$,  and to study the two-weight boundedness for the strong fractional maximal  operator

\begin{equation}\label{eq:StrongMaximal}
{\mathcal{M}}^{(S)}_{\alpha_1,\cdots \alpha_k}(\overrightarrow{f})(x) = \sup_{ \substack{Q^{(1)} \times \cdots \times Q^{(k)} \ni x\;\;\; }}
\frac{1}{ |Q^{(1)}|^{m-\alpha_1/n} \times \cdots \times |Q^{(k)}|^{m-\alpha_k/n} } \prod_{i=1}^m \int\limits_{Q^{(1)}\times\cdots \times  Q^{(k)}}  |f_i(y_i)|dy_i,
\end{equation}
where $0\leq \alpha_1,\cdots, \alpha_k < mn$, $ x= (x_1, \cdots, x_k) \in {\Bbb{R}}^{kn}$ ($k\geq 2$) and the supremum is taken over all products of cubes
$Q^{(1)}, \cdots, Q^{(k)} \in {\mathcal{Q}}$ containing $x_1, \cdots, x_k$ respectively.

Operator \eqref{eq:StrongMaximal} for $\alpha_1=\cdots \alpha_k=0$ and $n=1$ was introduced in \cite{GrLiPeTo}. In this case we have multi(sub)linear strong maximal operator denoted by ${\mathcal{M}^{(S)}}$ and defined with respect to rectangles in ${\Bbb{R}}^k$ with sides parallel to the coordinate
axes. In that paper the authors studied one and two--weight problems for ${\mathcal{M}}^{(S)}$.  In particular, they proved that the one-weight inequality holds if and only if the weight satisfies the strong $A_{\overrightarrow{p}}$ condition ($A_{\overrightarrow{p}}$ condition written with respect to rectangles in ${\Bbb{R}}^k$).

\begin{theorem}\label{mainresult-1}
Let $1<p_i<\infty$, $i=1, \cdots, m$. Suppose that $p < q <\infty$  and  $0< \alpha<mn$. Suppose also that  $w_i^{-p_i'}$ satisfies the reverse doubling condition on ${\Bbb{R}}^n$, $i=1, \cdots, m$. Then the inequality
\begin{equation}\label{eq:two-weightinequality}
\| u {\mathcal{M}}_{\alpha}(\overrightarrow{f})\|_{L^q({\Bbb{R}}^n)} \leq C \prod_{i=1}^{m} \|w_i f_i \|_{L^{p_i}({\Bbb{R}}^n)}
\end{equation}
holds for all $\overrightarrow{f} \in L^{p_1}_{w_1^{p_1}}({\Bbb{R}}^n) \times \cdots \times L^{p_m}_{w_m^{p_m}}({\Bbb{R}}^n)$ if and only if condition \eqref{eq:two-weightCondition} is satisfied.
\end{theorem}

\vskip+0.1cm
\begin{remark}
Since the reverse doubling condition is weaker than the $A_{\infty}$ condition (see e.g. \cite{GCRF}), Theorem \ref{mainresult-1} generalizes Corollary A for ${\mathcal{N}}_{\alpha}= {\mathcal{M}}_{\alpha}$.
\end{remark}

\begin{corollary}\label{corollary1}
Let conditions of Theorem \ref{mainresult-1} be satisfied and let, in addition, $u^q\in A_{\infty}({\Bbb{R}}^n)$. Then the inequality
$$ \| u {\mathcal{I}}_{\alpha}(\overrightarrow{f})\|_{L^q({\Bbb{R}}^n)} \leq C \prod_{i=1}^{m} \|w_i  f_i \|_{L^{p_i}({\Bbb{R}}^n)} $$
holds for all $\overrightarrow{f} \in L^{p_1}_{w_1^{p_1}}({\Bbb{R}}^n) \times \cdots \times L^{p_m}_{w_m^{p_m}}({\Bbb{R}}^n)$ if and only if
condition \eqref{eq:two-weightCondition} is satisfied.
\end{corollary}

\vskip+0.1cm

\begin{corollary}[\cite{Mo}] \label{corollary2}
Let $1<p_i<\infty$ ($i=1, \cdots, m$) and let  $p < q <\infty$. Suppose that  $0< \alpha<n$. We set $1/q= 1/p-\alpha/n$. Then the one-weight inequality
$$ \bigg(\int\limits_{{\Bbb{R}}^n}\Big[ \big| {\mathcal{N}}_{\alpha}(\overrightarrow{f})(x) \big| \Big( \prod_{i=1}^m w_i(x)\Big) \Big]^q dx \bigg)^{1/q} \leq C \prod_{i=1}^m \bigg( \int\limits_{{\Bbb{R}}^n} \Big[ \Big| f_i(x) \Big| w_i(x)\Big]^{p_i} dx \bigg)^{1/p_i},  $$
where ${\mathcal{N}}_{\alpha}$ is ${\mathcal{M}}_{\alpha}$ or ${\mathcal{I}}_{\alpha}$,
holds if and only if $\overrightarrow{w} \in A_{\overrightarrow{p}, q}$.
\end{corollary}

Regarding the strong  maximal operator we have:

 \vskip+0.2cm

\begin{theorem}\label{mainresult-2}
Let $1<p_i<\infty$ $(i=1, \cdots, m)$ and let  $p< q <\infty$. Assume that  $0< \alpha_1, \cdots, \alpha_k <mn$. Suppose that $u$ and $w_i$ are a.e. positive functions defined on ${\Bbb{R}}^{kn}$ such that $u^q$, $w_1^{p_1}, \cdots, w_m^{p_m}$ are weights. Suppose, in addition, that  $w_i$ are of the product type: $w_i= w_i^{(1)}\times \cdots \times w_i^{(k)}$, $i=1, \cdots, m$, where $w_i^{(1)}, \cdots, w_i^{(k)}$ are defined on ${\Bbb{R}}^n$. Assume also that $\big(w_i^{(1)}\big)^{-p'}, \cdots,  \big(w_i^{(k)}\big)^{-p'}\in RD({\Bbb{R}}^n)$.
Then the inequality
\vskip+0.2cm
$$ \| u {\mathcal{M}}^{(S)}_{\alpha_1, \cdots, \alpha_k}(\overrightarrow{f})\|_{L^q({\Bbb{R}}^{kn})} \leq C \prod_{i=1}^{m} \| w_i f_i \|_{L^{p_i} ({\Bbb{R}}^{kn})} $$
holds if and only if

\begin{equation} \label{eq:two-weightConditionStrong}
\sup_{Q_1,\cdots,  Q_k \in {\mathcal{Q}}}  |Q_1|^{\alpha_1/n-m}\times \cdots \times  |Q_k|^{\alpha_k/n-m} \bigg( \int\limits_{Q_1\times
\cdots \times Q_k} u^{q}(x) dx \bigg)^{1/q} \prod_{i=1}^m \bigg( \int\limits_{Q_1 \times \cdots \times Q_k} w_i^{-p'_i}(y_i) dy_i \bigg)^{1/p'_i}< \infty.
\end{equation}
\end{theorem}

Since for a constant weight function the reverse doubling condition is automatically satisfied, from Theorem \ref{mainresult-2} we have a characterization of the boundedness  for the operator ${\mathcal{M}}_{\alpha_1, \cdots, \alpha_k}$ from $L^p$ to $L^q_{u^q}$. In particular, we have

\begin{corollary}\label{solutionTraceProblem}
Let $p, q, \alpha_1, \cdots, \alpha_k$ satisfy the conditions of Theorem \ref{mainresult-2}. Then the inequality
$$\| u {\mathcal{M}}_{\alpha_1, \cdots,  \alpha_k}^{(S)}(\overrightarrow{f})\|_{L^q({\Bbb{R}}^{kn})} \leq c
\prod_{i=1}^{m} \|f_i \|_{L^{p_i}({\Bbb{R}}^{kn})}   $$
holds for all $\overrightarrow{f} \in L^{p_1}({\Bbb{R}}^{kn}) \times \cdots \times L^{p_m}({\Bbb{R}}^{kn})$
if and only if

\begin{equation} \label{eq:traceCondition}
\sup_{Q_1,\cdots, Q_k} \bigg( \int\limits_{Q_1\times  \cdots \times  Q_k}u^q (x) dx \bigg) |Q_1|^{q(\alpha_1/n-1/p)}\times \cdots \times  |Q_k|^{q(\alpha_k/n-1/p)} < \infty,
\end{equation}
where, as before, $\frac{1}{p}=  \sum_{i=1}^m \frac{1}{p_i}$.
\end{corollary}

Similar results were derived in \cite{KoMe} (see also \cite{KoMePe}, Ch.4) for sublinear maximal operators.

\vskip+0.1cm

It is natural to study the appropriate potential operator with product kernels:

$$ {\mathcal{I}}_{\alpha_1, \cdots, \alpha_k} (\overrightarrow{f})(x^{(1)}, \cdots, x^{(k)})= \int\limits_{ ({\Bbb{R}}^{kn})^m }
\frac{ \prod_{i=1}^m f_i (y^{(1)}_i,\cdots,  y^{(k)}_i) d y^{(1)}_1\cdots  d y^{(1)}_k\cdots d y^{(m)}_1,\cdots d y^{(m)}_k}   { \prod_{s=1}^k \big(\sum_{i=1}^m |x^{(s)}-y_i^{(s)}|\big)^{mn-\alpha_s} }, $$
where $x^{(i)}\in {\Bbb{R}}^n$, $i=1, \cdots, m$.
\vskip+0.2cm

It is easy to see that the following pointwise inequality holds
$$ {\mathcal{I}}_{\alpha_1,\cdots, \alpha_k} (\overrightarrow{f}) \geq c_n  {\mathcal{M}}^{(S)}_{\alpha_1,\cdots, \alpha_k}(\overrightarrow{f}), \;\;\; f_i\geq 0, \; i=1, \cdots, m.$$
\vskip+0.1cm

If $m=1$, the operator ${\mathcal{I}}_{\alpha_1, \cdots, \alpha_n}$ is the linear Riesz potential operator with product kernels denoted by $I_{\alpha_1, \cdots, \alpha_n}$. The one-weight criteria for this operator were established in \cite{Ko}, while the trace type inequality  was derived in \cite{KoMe} (see also \cite{KoMePe}, Ch.4).

\begin{theorem} \label{PotentialsProductKernels}
Let $p_i, q, \alpha_1, \cdots, \alpha_k$ satisfy the conditions of Theorem \ref{mainresult-2}. Suppose that $u^q$ satisfies $A_{\infty}({\Bbb{R}}^n)$ condition in each variable separately uniformly to another one. Then the trace type inequality
$$\|u {\mathcal{I}}_{\alpha_1, \cdots \alpha_k}(\overrightarrow{f})\|_{L^q({\Bbb{R}}^{kn})} \leq c
\prod_{i=1}^{m} \|w_if_i \|_{L^{p_i}({\Bbb{R}}^{kn})}   $$
holds if and only if \eqref{eq:two-weightConditionStrong} is true.
\end{theorem}

\begin{corollary} [One-weight inequality] \label{one-weight-strong}
Let $1<p_i<\infty$, $i=1, \cdots, m$. Suppose that $p < q <\infty$  and  $0< \alpha<mn$. We set $1/q=1/p-\alpha/n$. Suppose that $w_i$ are a.e. positive functions such that $w_i^{p_i}$ are weights on ${\Bbb{R}}^{kn}$ and $w_i= w_i^{(1)}\times \cdots \times w_i^{(k)}$, $i=1, \cdots, m$, where $w_i^{(j)}$ are defined on ${\Bbb{R}}^n$. Then the one-weight inequality
$$ \bigg(\int\limits_{ {\Bbb{R}}^{nk} }\Big[ \big| {\mathcal{S}}_{\alpha}(\overrightarrow{f})(x) \big| \Big( \prod_{i=1}^m w_i(x)\Big) \Big]^q dx \bigg)^{1/q} \leq C \prod_{i=1}^m \bigg( \int\limits_{{\Bbb{R}}^{nk}} \Big[ \Big| f_i(x) \Big| w_i(x)\Big]^{p_i} dx \bigg)^{1/p_i},  $$
where ${\mathcal{S}}_{\alpha}$ is ${\mathcal{M}}^{(S)}_{\alpha, \cdots, \alpha}$ or ${\mathcal{I}}_{\alpha, \cdots, \alpha}$,
holds if and only if $\overrightarrow{w} \in {\mathcal{A}}_{\overrightarrow{p}, q}({\Bbb{R}}^{kn})$, i.e.
\begin{equation}\label{one-weight-strong-condition}
\sup_{Q_1,\cdots, Q_k\in {\mathcal{Q}} } \bigg( \frac{1}{\prod_{i=1}^k|Q_i|} \int\limits_{Q_1 \times\cdots \times Q_k} \Big( \prod_{i=1}^m w_i\Big)^{q}\bigg)^{1/q} \prod_{i=1}^m \bigg( \frac{1}{\prod_{i=1}^k|Q_i|} \int\limits_{Q_1\times \cdots \times Q_k} w_i^{-p_i'}\bigg)^{1/p_i'}< \infty,
\end{equation}
where the supremum is taken over all cubes $Q_j \subset  {\Bbb{R}}^n$, $j=1, \cdots, k$.
\end{corollary}

Our next result concerns with the Fefferman-Stein type inequality for the (sub)linear strong fractional maximal operator ${\mathcal{M}}_{\alpha_1, \cdots, \alpha_k}^{(S)}$.

\begin{theorem} \label{Fefferman-Stein}
Let $p, q, \alpha_1, \cdots, \alpha_k$  satisfy the conditions of Theorem \ref{mainresult-2}. Suppose that $v$ is an a.e. positive function on ${\Bbb{R}}^{kn}$. Then there is a positive constant $C$ such that
 $$\| {\mathcal{M}}^{(S)}_{\alpha_1, \cdots, \alpha_k}(\overrightarrow{f})\|_{L_v^q({\Bbb{R}}^{kn})} \leq c
\prod_{i=1}^{m} \| (\overline{M}_{\alpha_1, \cdots, \alpha_k}v )^{p/(p_iqm)}f_i \|_{L^{p_i}({\Bbb{R}}^{kn})}   $$
holds, where $\overline{M}_{\alpha_1, \cdots, \alpha_m}$ is the sublinear strong fractional maximal operator defined on ${\Bbb{R}}^{kn}$  and given by
$$ \Big( \overline{M}_{\alpha_1, \cdots, \alpha_m}v\Big)( x^{(1)}, \cdots,   x^{(k)})$$
$$= \sup_{Q_1\ni x^{(1)}, \cdots,  Q_k \ni x^{(k)}} |Q_1|^{q\big( \frac{\alpha_1}{n}-\frac{1}{p}\big)}\times \cdots \times   |Q_k|^{q\big( \frac{\alpha_k}{n}-\frac{1}{p}\big)} \int_{Q_1\times, \cdots, \times  Q_k} v(y^{(1)},\cdots,  y^{(k)}) dy^{(1)}\cdots  dy^{(k)}. $$
\end{theorem}

\section{Proofs of the Main Results}

In this section we give proofs of the main results of this paper. First we formulate the following auxiliary statement:
\vskip+0.2cm
\begin{lemma}\cite{Mo}
Let ${\mathcal{M}}_{\alpha}^k$ be defined as ${\mathcal{M}}_{\alpha}$ but over cubes with side length less than or equal to $2^k$. Suppose that  $Q_k:= [-2^{k+2}, 2^{k+2})^n$, $\tau_tg(x):= g(x-t)$ and $\overrightarrow{\tau}_t (\overrightarrow{f}):=
(\tau_tf_1, \cdots, \tau_tf_m).$  Then there is a positive constant $c$ depending only on $n$, $m$, $\alpha$ and $q$ such that for each $k\in {\Bbb{Z}}$, $\overrightarrow{f}$ and every  $x\in {\Bbb{R}}^n$,
\begin{equation} \label{eq:Fefferman-SteinPointwise}
\big( {\mathcal{M}}^k_{\alpha} (\overrightarrow{f})(x) \big)^q \leq \frac{c}{|Q_k|} \int\limits_{Q_k} \Big[ \Big( \tau_{-t}\circ {\mathcal{M}}_{\alpha}^{(d)} \circ \overrightarrow{\tau}_t\Big)\big( \overrightarrow{f}\big)(x)\Big]^q dt
\end{equation}
\end{lemma}

{\em Proof of Theorem \ref{mainresult-1} }. {\em Necessity.} Let us first show that
\begin{equation}\label{eq:weightisfinite}
\int\limits_{Q}w_{i}^{-p_i'}(y_i) dy_i<\infty,
\end{equation}
for all  $i=1, \cdots, m$ and all cubes $Q$. Indeed, suppose that

\begin{equation}\label{eq:weightisinfinite}
\int\limits_{Q}w_{i_1}^{-p_{i_1}'}(y_{i_1}) dy_{i_1}= \cdots = \int\limits_{Q}w_{i_k}^{-p_{i_k}'}(y_{i_k}) dy_{i_k}= \infty,
\end{equation}
for some $1\leq i_1< \cdots < i_k\leq m$. Then by the duality arguments, we have that there are  non--negative functions $g_{i_1}, \cdots g_{i_k}$  such that
$$ \int\limits_{Q}w_{i_j}^{-1} g_{i_j} =\infty, \;\;\;\; j=1, \cdots, k,  $$
and $g\in L^{p_{i_j}}(Q)$,  $j=1, \cdots, k$.  Let us assume that $ \overrightarrow{f}= (f_1, \cdots, f_m), $
where $f_{j} =w^{-1}_{i_j}g_{i_j}\chi_{Q} $, $j=1, \cdots, k$  and $f_n= \chi_{Q}$ for other $n$.  Then

$$ \| u {\mathcal{M}}_{\alpha} (\overrightarrow{f}) \|_{L^q({\Bbb{R}}^n)} \geq \frac{1}{|Q|^{m-\alpha/n}}
\bigg( \prod_{j=1}^k \int\limits_{Q} w_{i_j}^{-1} g_{i_j}\bigg) |Q|^{m-1} \|u\|_{L^q(Q)}  = \infty. $$
On the other hand,
$$ \prod_{i=1}^{m} \|f_i w_i \|_{L^{p_i}({\Bbb{R}})} < \infty $$
because $w_i^{p_i}$ are locally integrals and $g_{i_j} \in L^{p_{i_j}}({\Bbb{R}}^n)$, $j=1, \cdots, k$.

{\em Sufficiency.} First we show that the two-weight inequality
 \begin{equation}\label{eq:dyadictwoweight}
 \| u {\mathcal{M}}_{\alpha}^{(d)} (\overrightarrow{f})\|_{L^q({\Bbb{R}}^n)} \leq C \prod_{i=1}^{m} \| w_i f_i \|_{L^{p_i}({\Bbb{R}}^n)}
 \end{equation}
holds if \eqref{eq:two-weightCondition} holds provided that $w_i^{-p'_i} \in RD^{(d)}({\Bbb{R}}^n)$, $i=1, \cdots, m$, where ${\mathcal{M}}_{\alpha}^{(d)}$ is a dyadic multi(sub)linear fractional maximal operator defined  by


$$ ({\mathcal{M}}_{\alpha}^{(d)}\overrightarrow{f})(x) =\sup_{Q\ni {\mathcal{Q}}, Q\ni x} \prod_{i=1}^m \frac{1}{ |Q|^{1-\alpha/(nm)}} \int\limits_{Q} |f_i(y_i)| d y_i, \;\;\; 0< \alpha< mn.  $$

To show that \eqref{eq:two-weightCondition} implies \eqref{eq:dyadictwoweight} we argue as follows:  for every $x\in {\Bbb{R}}^n$ let us take $Q_x \in {\mathcal{D}}$ such that $Q_x\ni x$ and
\begin{equation}\label{eq:pointwiseestimate}
\big( {\mathcal{M}}^{(d)}_{\alpha}\overrightarrow{f}\big)(x) \leq \frac{2}{|Q_x|^{m-\alpha/n}} \prod_{i=1}^m\int\limits_{Q_x}
|f_i(y_i)| dy_i.
\end{equation}

Without loss of generality we can assume, for example,  that $f_i$, $i=1, \cdots, m$  are non--negative, bounded and have compact supports.

Let us introduce a set
$$ F_{Q}= \{ x\in {\Bbb{R}}^n: x\in Q\;  \text{and} \;\; \eqref{eq:pointwiseestimate}\;\;  holds\;\;   for\;\;  Q\}. $$

It is obvious that $F_Q\subset Q$ and ${\Bbb{R}}^n =\cup_{Q\in {\mathcal{D}}} F_Q$.

Further, we have

$$ \int\limits_{ {\Bbb{R}}^n} \big( {\mathcal{M}}^{(d)}_{\alpha}\overrightarrow{f}\big)(x)\big)^q(x) u^q(x) dx \leq  \sum_{Q\in {\mathcal{D}}} \int\limits_{F_Q} \big({\mathcal{M}}^{(d)}_{\alpha}\overrightarrow{f}\big)^q u^q(x) dx $$

$$ \leq  2^q \sum_{Q\in {\mathcal{D}}} |Q|^{(\alpha/n-m)q} \bigg( \int\limits_{Q} u^q(x) dx \bigg) \bigg( \prod_{i=1}^m \int\limits_{Q} f_i(y_i) dy_i \bigg)^q $$

$$ \leq c \sum_{Q\in {\mathcal{D}}} \bigg( \prod_{i=1}^m \int\limits_{Q} w_i^{-p'}(y_i) dy_i\bigg)^{-q/p_i'} \bigg( \prod_{i=1}^m \int\limits_{Q} f_i(y_i) dy_i \bigg)^q$$
(using  H\"older's inequality  $\sum_{k}a^{(1)}_k\times \cdots \times   a^{(m)}_k \leq  \prod_{j=1}^m \big( \sum_k (a^{(j)}_k)^{p_j/p}\big)^{p/p_j}$ for positive sequences $\{ a^{(j)}_k\}$,  $j=1,\cdots, m$)

$$ \leq \bigg[\sum_{Q\in {\mathcal{D}}} \bigg( \int\limits_{Q} w_1^{-p_1'}(y_1) dy_1\bigg)^{-(q p_1)/(p p_1')} \bigg( \int\limits_{Q} f_1(y_1) dy_1\bigg)^{q p_1/p}\bigg)^{p/p_1} $$
$$ \times \cdots \times \bigg( \sum_{Q\in {\mathcal{D}}} \bigg( \int\limits_{Q} w_m^{-p_m'}(y_m) dy_m\bigg)^{-(qp_m)/(p p_m')} \bigg( \int\limits_{Q} f_m(y_m) dy_m\bigg)^{q p_1/p}\bigg]^{p/p_1}   $$
(Applying Theorem G for the exponents $(p_i, q p_i/p)$, $i= 1, \cdots, m$)
$$\leq c \|f_1 w_1\|^q_{L^{p_1}({\Bbb{R}}^n)} \times \cdots \times \|f_m w_m\|^q_{L^{p_m}({\Bbb{R}}^n)}. $$

Now we pass from ${\mathcal{M}}^{(d)}_{\alpha}$ to ${\mathcal{M}}_{\alpha}$. By using \eqref{eq:Fefferman-SteinPointwise} we see that
$$ \| u {\mathcal{M}}_{\alpha}(\overrightarrow{f})\|_{L^q({\Bbb{R}}^n)} \leq \sup_{t} \| u \Big( \tau_{-t} \circ {\mathcal{M}}_{\alpha}^{(d)} \circ \overrightarrow{\tau}_{t}\Big) \big(\overrightarrow{f}\big) \|_{L^q({\Bbb{R}}^n)}.$$

Observe now that  a pair $(\tau_t u, \overrightarrow{\tau}_t \overrightarrow{w})$ satisfy condition \eqref{eq:two-weightCondition} independent of $t$. Notice also that the weight functions $\tau_t \circ w_i $ satisfy the reverse doubling condition. Hence,

$$   \| u \Big( \tau_{-t} \circ {\mathcal{M}}_{\alpha}^{(d)} \circ \overrightarrow{\tau}_{t} \Big)\big(\overrightarrow{f}\big) \|_{L^q({\Bbb{R}}^n)} \leq c \prod_{i=1}^m \|w_if_i\|_{L^{p_i}({\Bbb{R}}^n)} $$
with the constant $c$ independent of $t$. Finally, we conclude that
$$ \int\limits_{{\Bbb{R}}^n} \big( {\mathcal{M}}_{\alpha} \overrightarrow{f}\big)(x)\big)^q(x) u^q(x) dx  \leq \sup_{t} \Big\| u \Big( \tau_{-t} \circ {\mathcal{M}}_{\alpha}^{(d)} \circ \overrightarrow{\tau}_t\Big) \big( \overrightarrow{f}\big) \Big\|_{L^q({\Bbb{R}}^n)} \leq c \prod_{i=1}^m \|f_i w_i\|_{L^{p_i}({\Bbb{R}}^n)}. $$

$\Box$

\vskip+0.2cm

\begin{lemma}\label{A-infty}
Let $0<q<\infty$ and let $0<\alpha< mn$. Suppose that a weight function $v$ satisfies $A_{\infty}({\Bbb{R}}^n)$. Then there is a positive constant $C$ such that for all $\overrightarrow{f}$ the inequality
$$\int\limits_{{\Bbb{R}}^n} |{\mathcal{I}}_{\alpha}(\overrightarrow{f})(x)|^q v(x) dx \leq c
\int\limits_{{\Bbb{R}}^n} \big( {\mathcal{M}}_{\alpha}(\overrightarrow{f})(x)\big)^q v(x) dx$$
holds.
\end{lemma}
Corollary \ref{corollary1} follows immediately from Theorem \ref{mainresult-1} and Lemma \ref{A-infty}.
\vskip+0.2cm

Corollary \ref{corollary2} for ${\mathcal{M}}_{\alpha}$ is a direct consequence of Lemma \ref{MuckenhouptProperties} and the fact that the Muckenhoupt condition implies the reverse doubling condition (see e.g. \cite{GCRF}, \cite{Wh}). Now the result for the operator ${\mathcal{I}}_{\alpha}$ follows by applying Lemma \ref{A-infty}
\vskip+0.2cm

To prove Theorem 3.4 we introduce a strong dyadic fractional operator

$$
\big( {\mathcal{M}}_{\alpha, \beta}^{(S),(d)}(\overrightarrow{f})(x) = \sup_{\substack{ Q^{(1)}\times Q^{(2)}\ni x,\;} {Q^{(1)}, Q^{(2)} \in {\mathcal{D}} }}\;\;
\frac{1}{ |Q^{(1)}|^{m-\alpha/n} |Q^{(2)}|^{m-\beta/n} } \prod_{i=1}^m \int\limits_{Q^{(1)}\times Q^{(2}}  |f_i(y_i)|dy_i. $$
We need also the following notation

\begin{equation}\label{eq:StrongMaximal-truncated}
\big( {\mathcal{M}}_{\alpha, \beta}^{(S),(k)}(\overrightarrow{f})(x) = \sup_{\substack{Q^{(1)}\times Q^{(2)} \ni x,} { |Q^{(1)}|, |Q^{(2)}|\leq  2^k}}
\frac{1}{ |Q^{(1)}|^{m-\alpha/n} |Q^{(2)}|^{m-\beta/n} } \prod_{i=1}^m \int\limits_R  |f_i(y_i)|dy_i;
\end{equation}
$$ Q_k:= [-2^{k+2}, 2^{k+2})^{2n}, \;\; \tau_{t,\delta} g(x^{(1)}, x^{(2)}):= g(x^{(1)}-t, x^{(2)}-\delta)$$

$$ \overrightarrow{\tau}_{t, \delta} (\overrightarrow{f}):= (\tau_{t, \delta} f_1, \cdots, \tau_{t, \delta} f_m),$$
where $g$ and $f$ are measurable functions on ${\Bbb{R}}^{2n}$.

\begin{lemma}\label{strong-fefferman-stein}
There exists a positive constant $c$ depending only on $n$, $m$, $\alpha$, $\beta$ and $q$ such that for all $k\in {\Bbb{Z}}$, $\overrightarrow{f}$ and every  $x\in {\Bbb{R}}^{2n}$,
\begin{equation} \label{eq:StrongFefferman-SteinPointwise}
\big( {\mathcal{M}}^{(S),(k)}_{\alpha, \beta} (\overrightarrow{f})(x) \big)^q \leq \frac{c}{|Q_k|^2} \int\limits_{Q_k^2} \bigg[ \Big( \tau_{-t, -\delta}\circ {\mathcal{M}}_{\alpha, \beta}^{(S),(d)} \circ \overrightarrow{\tau}_{t, \delta} \Big)\big( \overrightarrow{f}\big)(x)\bigg]^q dt d\delta.
\end{equation}
\end{lemma}
{\em Proof} of this lemma follows in the same way as in the case of (sub)linear case (see \cite{KoMe}, \cite{KoMePe}, P. ); therefore we omit the details.

{\em Proof of Theorem \ref{mainresult-2}.} For simplicity we prove the theorem for $k=2$. The proof for other $k$ is similar.

{\em Necessity} follows in the same way as in the case of the operator ${\mathcal{M}}_{\alpha}$.

To show {\em Sufficiency} we assume that $f_i$ are non-negative and  bounded with compact supports.

For every $(x^{(1)}, x^{(2)})\in {\Bbb{R}}^{2n}$, we take  $Q^{(1)}, Q^{(2)} \in {\mathcal{D}}$ such that $Q^{(1)}\ni x^{(1)}$, $Q^{(2)}\ni x^{(2)}$ and

\begin{equation}\label{eq:strongpointwiseestimate}
{\mathcal{M}}^{(S),(d)}_{\alpha_1, \alpha_2} \big( \overrightarrow{f}\big) (x^{(1)}, x^{(2)}) \leq \frac{2}{ |Q^{(1)}|^{m-\alpha_1/n} |Q^{(1)}|^{m-\alpha_2/n} } \prod_{i=1}^m \int\limits_{ Q^{(1)}\times Q^{(2)} }
f_i(y^{(1)}_i, y^{(2)}_i) dy^{(1)}_i dy^{(2)}_i.
\end{equation}

We construct the set:

$$ F_{Q^{(1)}, Q^{(2)}}:= \{ x= (x^{(1)}, x^{(2)}) \in {\Bbb{R}}^{2n}: x^{(1)}\in Q^{(1)}, x^{(2)}\in Q^{(2)} \text{and} \; \eqref{eq:strongpointwiseestimate}\;\; \text{holds}\;\;  \text{for}\;\;  Q^{(1)}\;\; \text{and}\;\; Q^{(2)}\}. $$

Since $F_{Q^{(1)}, Q^{(2)}} \subset Q^{(1)}\times Q^{(2)}$, arguing as in the proof of Theorem \ref{mainresult-1} we find that

$$ \int\limits_{ {\Bbb{R}}^{2n} } \big( {\mathcal{M}}^{(d)}_{\alpha_1, \alpha_2} (\overrightarrow{f}) \big)(x^{(1)}, x^{(2)}) \big)^q u^q ( x^{(1)}, x^{(2)}) dx^{(1)} x^{(2)} $$

$$\leq  \sum_{Q_1, Q_2 \in {\mathcal{D}} } \int\limits_{F_{Q^{(1)}, Q^{(2)}}} \big( {\mathcal{M}}^{(d)}_{\alpha_1, \alpha_2}(\overrightarrow{f}) (x^{(1)}, x^{(2)}) \big)^q  u^q (x^{(1)}, x^{(2)}) dx^{(1)} x^{(2)} $$

$$ \leq  2^q \sum_{Q^{(1)}, Q^{(2)} \in {\mathcal{D}} } |Q|^{(\alpha/n-m)q} \bigg( \int\limits_{Q} u^q(x^{(1)}, x^{(2)}) dx^{(1)}x^{(2)} \bigg) \bigg( \prod_{i=1}^m \int\limits_{Q^{(1)}\times Q^{(2)}} f_i(y^{(1)}_i, y^{(2)}_i) dy^{(1)}_i dy^{(2)}_i \bigg)^q $$

$$ \leq c \sum_{Q^{(1)}, Q^{(2)} \in {\mathcal{D}}} \prod_{i=1}^m \bigg(\int\limits_{ Q^{(1)} \times Q^{(2)}} w_i^{-p'}(y_i) dy_i\bigg)^{-q/p_i'} \bigg( \prod_{i=1}^m  \int\limits_{Q} f_i(y^{(1)}_i, y^{(2)}_i) dy^{(1)}_i dy^{(2)}_i\bigg)^q $$

(due to  H\"older's inequality  $\sum_{k}a^{(1)}_k\times \cdots \times   a^{(m)}_k \leq  \prod_{j=1}^m \big( \sum_k (a^{(j)}_k)^{p_j/p}\big)^{p/p_j}$ for positive sequences $\{ a^{(j)}_k\}$,  $j=1,\cdots, m$)

$$ \leq \bigg[\sum_{ Q^{(1)}, Q^{(2)} \in {\mathcal{D}} } \bigg( \int\limits_{ Q^{(1)}} \big(w^{(1)}_1\big)^{-p_1'}(y^{(1)}_1) dy^{(1)}_1 \bigg)^{-(q p_1)/(p p_1')} \bigg( \int\limits_{Q^{(2)}} \big( w^{(2)}_1\big)^{-p_1'} (y^{(2)}_1) dy^{(2)}_1 \bigg)^{-(qp_1)/(pp_1')} $$
$$\times \bigg( \int\limits_{Q^{(1)}} \bigg( \int\limits_{Q^{(2)}} f_1(y^{(1)}_1, y^{(2)}_1) dy^{(1)}_1 \bigg) dy^{(2)}_1 \bigg)^{q p/ p_1}\bigg]^{p_1/p} $$

$$ \times \cdots \times \bigg[ \sum_{Q^{(1)}, Q^{(2)} \in {\mathcal{D}} } \bigg( \int\limits_{Q^{(1)}} \big( w^{(1)}_m\big)^{-p_m'}(y^{(1)}_m) dy^{(1)}_m \bigg)^{-(q p_m)/(p p_m')} \bigg( \int\limits_{Q^{(2)}}\big(w^{(2)}_m\big)^{-p_1'}(y^{(2)}_m) dy^{(2)}_m \bigg)^{-(q p_1)/(p p_m')} $$
$$ \times \bigg( \int\limits_{Q^{(1)}} \bigg( \int\limits_{Q^{(2)}} f_1(y^{(1)}_m, y^{(2)}_m) dy^{(1)}_m\bigg) \; dy^{(2)}_m \bigg)^{qp/p_m}\bigg]^{p_m/p} $$
(by virtue of Theorem  G for the exponents $(p_i, q p_i/p)$, $i= 1, \cdots, m$, with respect to the first variable)

$$  \leq \bigg[ \sum_{Q^{(2)} \in {\mathcal{D}}} \bigg( \int\limits_{Q^{(2)}} \big(w^{(2)}_1\big)^{-p_1'}(y^{(2)}_1)  dy_1^{(2)} \bigg)^{-(qp_1)/(pp_1')} \bigg( \int\limits_{{\Bbb{R}}^n} \big(w^{(1)}_1\big)^{p_1}(y^{(1)}_1) \bigg( \int\limits_{Q^{(2)}} f_1(y^{(1)}_1, y^{(2)}_1) dy^{(2)}_1\bigg)^{p_1} dy^{(1)}_1 \bigg)^{q/p} \bigg]^{p/p_1}$$

$$ \times \cdots \times \bigg[\sum_{Q^{(2)} \in {\mathcal{D}}} \bigg( \int\limits_{Q^{(2)}} \big(w^{(2)}_m\big)^{-p_m'}(y^{(2)}_m) dy^{(2)}_m \bigg)^{-(qp_m /pp_m')} \bigg( \int\limits_{{\Bbb{R}}^n} \big( w^{(1)}_m\big)^{p_m}(y^{(1)}_m) \bigg( \int\limits_{Q^{(2)}} f_m(y^{(1)}_m, y^{(2)}_m) dy^{(2)}_1\bigg)^{p_m} dy^{(1)}_m \bigg)^{q/p}\bigg]^{p/p_m} $$
(using generalized Minkowski's inequality)

$$\leq  \bigg[\sum_{Q^{(2)} \in {\mathcal{D}} }\bigg( \int\limits_{Q^{(2)}} \big(w^{(2)}_1\big)^{-p_1'}(y^{(2)}_1) dy_1^{(2)} \bigg)^{-(qp_1/pp_1')} \bigg( \int\limits_{Q^{(2)}} \bigg( \int\limits_{{\Bbb{R}}^n} f_1^{p_{1}} ( y^{(1)}_1, y^{(2)}_1) \big(w^{(1)}_1\big)^{p_1}(y^{(1)}_1) dy^{(1)}_1\bigg)^{1/p_1} dy^{(2)}_1 \bigg)^{qp_1/p}\bigg]^{p/p_1} $$

$$ \times \cdots \times \bigg[\sum_{Q^{(2)} \in {\mathcal{D}}} \bigg( \int\limits_{Q^{(2)}} \big(w^{(2)}_m\big)^{-p_m'}(y^{(2)}_m) dy_m^{(2)} \bigg)^{-(qp_m/pp_m')} $$

$$ \times \bigg( \int\limits_{Q^{(2)}} \bigg( \int\limits_{{\Bbb{R}}^n} f_m^{p_{m}} ( y^{(1)}_m, y^{(2)}_m) \big(w^{(1)}_m\big)^{p_m}(y^{(1)}_m) dy^{(1)}_m\bigg)^{1/p_m} dy^{(2)}_m \bigg)^{qp_m/p}\bigg]^{p/p_m} $$
(by virtue of Theorem G for the pairs of exponents: $(p_i, q p_i/p)$, $i= 1, \cdots, m$  with respect to the second variable)

$$\leq c \|f_1 w_1\|^q_{L^{p_1}({\Bbb{R}}^{2n})} \times \cdots \times \|f_m w_m\|^q_{L^{p_m}({\Bbb{R}}^{2n})}. $$

Arguing now as in the proof of Theorem \ref{mainresult-1} and using Lemma \ref{strong-fefferman-stein} we can pass from the dyadic strong maximal function ${\mathcal{M}}_{\alpha_1, \alpha_2}^{(S)(d)}$ to the
operator ${\mathcal{M}}^{(S)}_{\alpha_1, \alpha_2}$ and get the desired result.

$\Box$

Theorem \ref{PotentialsProductKernels} follows from Theorem \ref{mainresult-2} and by using Lemma \ref{A-infty} with respect to each variable separately.

{\em Proof of Corollary \ref{one-weight-strong}}. We assume that  $k=2$ because the proof for higher $k$ is similar. Suppose that $w_i= w^{(1)}_i w^{(2)}_i$, where  $ w^{(1)}_i$ and $ w^{(2)}_i$ are defined on  ${\Bbb{R}}^n$. Assume that ${\mathcal{S}}_{\alpha_1, \alpha_2}= {\mathcal{M}}_{\alpha_1, \alpha_2}$. First observe that condition \eqref{one-weight-strong-condition} implies that

$$
\sup_{Q_1 \in {\mathcal{Q}}} \bigg( \frac{1}{|Q_1|} \int\limits_{Q_1} \Big[\prod_{i=1}^m w^{(1)}_i (x^{(1)}) \Big]^{q}dx^{(1)}\bigg)^{1/q} \prod_{i=1}^m \bigg( \frac{1}{|Q_1|} \int\limits_{Q_1}\big( w^{(1)}_{i}\big)^{-p'_i}(x^{(1)}) dx^{(1)} \bigg)^{1/p'}< \infty, $$

and
$$
\sup_{Q_2\in {\mathcal{Q}}} \bigg( \frac{1}{|Q_2|} \int\limits_{Q_2} \Big[ \prod_{i=1}^m w^{(2)}_i(x^{(2)})\Big]^{q} dx^{(2)} \bigg)^{1/q} \prod_{i=1}^m \bigg( \frac{1}{|Q_2|}\int\limits_{Q_2} \big( w^{(2)}_{i}\big)^{-p'_i}(x^{(2)}) dx^{(2)}\bigg)^{1/p_i'}< \infty. $$

Consequently, by Lemma \ref{MuckenhouptProperties} we have that $\big(w^{(1)}_i\big)^{-p_i'}, \big( w^{(2)}_i \big)^{-p_i'} \in A_{mp'_i}({\Bbb{R}}^n)$. Hence, (see e.g., \cite{GCRF} or \cite{Wh}) $\big(w^{(1)}_i\big)^{-p_i'}, \big( w^{(2)}_i \big)^{-p_i'}$  satisfy the reverse doubling condition on ${\Bbb{R}}^n$. Now the result follows from
Theorem \ref{mainresult-2}.

The  proof for ${\mathcal{S}}_{\alpha_1, \alpha_2}= {\mathcal{I}}_{\alpha_1, \alpha_2}$, follows from Lemmas  \ref{MuckenhouptProperties} and \ref{A-infty}.

$\Box$

{\em Proof of }Theorem \ref{Fefferman-Stein}. As before it is enough to prove the inequality

$$   \|{\mathcal{M}}^{(S), (d)}_{\alpha, \beta}(\overrightarrow{f})\|_{L^q_{v}({\Bbb{R}}^{2n})} \leq c
\prod_{i=1}^{m} \| \big(\overline{M}_{\alpha, \beta}v \big)^{p/(p_iqm)} f_i \|_{L^{p_i}({\Bbb{R}}^{2n})},  $$
where ${\mathcal{M}}^{(S), (d)}_{\alpha_1, \alpha_2}$ is a dyadic strong fractional (sub)linear maximal operator defined above. Constructing the sets $F_{Q^{(1)}, Q^{(2)}}$ for dyadic cubes $Q^{(1)}$ and $Q^{(2)}$ (see the proof of Theorem \ref{mainresult-2}) we have

$$ \int\limits_{{\Bbb{R}}^{2n}} \Big( {\mathcal{M}}^{(S), (d)}_{\alpha_1, \alpha_2} (\overrightarrow{f} )(x^{(1)}, x^{(2)}) \Big)^{q} v(x^{(1)}, x^{(2)}) dx^{(1)} dx^{(2)} \leq 2^q \sum_{Q^{(1)}, Q^{(2)}\in {\mathcal{D}}}\bigg( \int_{Q^{(1)} \times  Q^{(2)}} v(x^{(1)}, x^{(2)}) dx^{(1)} dx^{(2)}\bigg) $$

$$\times |Q^{(1)}|^{(\alpha_1/n-m)q} |Q^{(2)}|^{(\alpha_2/n-m)q} \prod_{i=1}^m \bigg( \int\limits_{Q^{(1)} \times Q^{(2)}} |f_i(y^{(1)}_i,  y^{(2)}_i)| dy^{(1)}_i  dy^{(2)}_i)\bigg)^q  $$

$$ =2^q \sum_{Q^{(1)}, Q^{(2)}\in {\mathcal{D}}}  |Q^{(1)}|^{(\alpha_1/n-m)q} |Q^{(2)}|^{(\alpha_2/n-m)q} $$

$$ \times \prod_{i=1}^m \bigg( \int\limits_{Q^{(1)}\times Q^{(2)}} |f_i(y^{(1)}_i,  y^{(2)}_i)| \bigg( \int_{Q^{(1)}\times  Q^{(2)}} v(x^{(1)}, x^{(2)}) dx^{(1)} dx^{(2)}\bigg)^{1/(qm)} dy^{(1)}_i  dy^{(2)}_i)\bigg)^q $$

$$ \leq c \sum_{Q^{(1)}, Q^{(2)} \in {\mathcal{D}}}  |Q^{(1)}|^{q(1/p-m)} |Q^{(2)}|^{(1/p-m)} \prod_{i=1}^m\bigg( \int\limits_{Q^{(1)}\times Q^{(2)}} |f_i(y^{(1)}_i,  y^{(2)}_i)| \big(  {\overline{M}}_{\alpha_1, \alpha_2} (v) (y^{(1)}_i, y^{(2)}_i) \big)^{1/(qm)}  dy^{(1)}_i  dy^{(2)}_i)\bigg)^q $$
(using the fact that $|Q^{(1)}\times Q^{(2)}|^{q(1/p-m)}= \prod_{i=1}^m  |Q^{(1)}\times Q^{(2)}|^{-q/p'_i}$)

$$  = c \sum_{Q^{(1)}, Q^{(2)}\in {\mathcal{D}}}  \prod_{i=1}^m |Q^{(1)}|^{-q/p'_i} |Q^{(2)}|^{-q/p'_i} $$

$$ \times \prod_{i=1}^m\bigg( \int\limits_{Q^{(1)} \times Q^{(2)}} |f_i(y^{(1)}_i,  y^{(2)}_i)| \big(  {\overline{M}}_{\alpha_1, \alpha_2} (v) (y^{(1)}_i, y^{(2)}_i) \big)^{1/(qm)}  dy^{(1)}_i  dy^{(2)}_i)\bigg)^q$$
(by virtue  of  H\"older's inequality)

$$ \leq c \bigg[\sum_{Q^{(1)}, Q^{(2)} \in {\mathcal{D}}}   |Q^{(1)}|^{-(qp_1/pp'_1)} |Q^{(2)}|^{-(qp_1)/(pp'_1)} \bigg( \int\limits_{Q^{(1)} \times Q^{(2)}} |f_1(y^{(1)}_1,  y^{(2)}_1)| \big(  {\overline{M}}_{\alpha_1, \alpha_2} (v) (y^{(1)}_1, y^{(2)}_1) \big)^{1/(qm)}  dy^{(1)}_1  dy^{(2)}_1)\bigg)^{qp_1/p}\bigg]^{p/p_1}$$

$$ \times \cdots \times \bigg[\sum_{Q^{(1)}, Q^{(2)}\in {\mathcal{D}} }   |Q^{(1)}|^{-(qp_m/pp'_m)} |Q^{(2)}|^{-(qp_m/pp'_m)} $$

$$ \times \bigg( \int\limits_{Q^{(1)} \times Q^{(2)}} |f_m(y^{(1)}_m,  y^{(2)}_m)| \big(  {\overline{M}}_{\alpha_1, \alpha_2} (v) (y^{(1)}_m, y^{(2)}_m) \big)^{1/(qm)}  dy^{(1)}_m  dy^{(2)}_m) \bigg)^{qp_m/p}\bigg]^{p/p_m}$$

(by applying twice Theorem G for $w\equiv 1$  with respect to each variable separately and Minkowski's generalized inequality)

$$ \leq c \prod_{i=1}^{m} \Big\| \Big( \overline{M}_{\alpha_1, \alpha_2} v\Big)^{p/(p_i qm)} f_i \Big\|_{L^{p_i}({\Bbb{R}}^{2n})}. $$

$\Box$
\vskip+0.2cm

\begin{remark}
It is easy to see that Corollary \ref{solutionTraceProblem} can be also obtained from Theorem \ref{Fefferman-Stein}.
\end{remark}

{\bf Acknowledgement.} The first and third authors were partially supported by the Shota Rustaveli National Science Foundation Grant (Contract Numbers D/13-23 and 31/47).


\end{document}